\documentclass[11pt]{article}
\usepackage{amsmath,amsthm}

\def\refeq#1{\if\workingver y(\ref{#1})-[[#1]]\else(\ref{#1})\fi}
\def\refth#1{\if\workingver y\ref{#1}-[[#1]]\else\ref{#1}\fi}
\def\mylabel#1{\if\workingver y\label{#1}{\bf\ \ [[#1]]\ \ }
\else\label{#1}\fi}
\def\mybibitem#1{\if\workingver y\bibitem{#1}{\bf\ \ [[#1]]\ \ }
\else\bibitem{#1}\fi}

\newtheorem{thm}{Theorem}

\newtheorem{algo}[thm]{Algorithm}

\newtheorem{conj}[thm]{Conjecture}

\def\begeq#1{\begin{equation}\mylabel{#1}}
\def\endeq{\end{equation}}

\def\begalg{\begin{alg}}
\def\endalg{\end{alg}}

\newcounter{li}

\def\begalg#1{\begin{algo}\mylabel{#1}\normalshape:\\}

\let\workingver=n

\font\tenrm=cmr10

\begin{document}
\title{\bf On the prime power factorization of $n!$}
\author{Florian Luca $^1$ \hskip15pt {\rm and}\hskip15pt
Pantelimon St\u anic\u a $^2$\thanks{Also associated with the Institute of
Mathematics of Romanian Academy, Bucharest, Romania.}\\
$^1$ \tenrm Instituto de Matematicas de la UNAM, Campus Morelia\\
 \tenrm Apartado Postal 61-3 (Xangari) CP 58 089, Morelia, Michoacan, Mexico\\
$^2$ \tenrm Department of Mathematics, Auburn University Montgomery\\
\tenrm Montgomery, AL 36124-4023, USA; e-mail: pstanica@mail.aum.edu}

\date{March 25, 2003
}
\maketitle

%
%

\begin{abstract}
In this paper we prove two results.
The first theorem uses a paper of Kim \cite{K} to
show that for fixed primes $p_1,\ldots,p_k$,
and for fixed integers $m_1,\ldots,m_k$, with $p_i\not|m_i$,
the numbers $(e_{p_1}(n),\ldots,e_{p_k}(n))$ are uniformly distributed modulo
$(m_1,\ldots,m_k)$, where $e_p(n)$ is the order of the prime $p$ in the
factorization of $n!$. That implies one of Sander's conjecture from \cite{S},
for any set of odd primes. Berend \cite{B} asks to find the fastest growing
function $f(x)$ so that for large $x$ and any given finite sequence
$\varepsilon_i\in \{0,1\},\ i\le f(x)$, there exists
$n<x$ such that the congruences $e_{p_i}(n)\equiv \varepsilon_i\pmod 2$
hold for all $i\le f(x)$. Here, $p_i$ is the $i$th prime number.
 In our second result, we are able to show that $f(x)$ can be taken
to be at least $c_1 (\log x/(\log\log x)^6)^{1/9}$,
 with some absolute constant $c_1$, provided that only the first odd
prime numbers are involved.
\end{abstract}
%
%
%

\section{Introduction}
\label{Introduction}


\par\noindent For a prime number $p$ and
every positive integer $m$ let $e_p(n)$
be the power at which the prime number $p$ appears in the prime
factorization of $n!$. In \cite{EG}, it was asked whether for every
fixed positive integer $k$ there exists a positive integer $n$ so that
all the numbers $e_{p_i}(n)$ are even, where $p_1<p_2<\dots<p_k$
denote the first $k$ prime numbers. The above question was answered
in the affirmative by Berend in \cite{B}. In fact, Berend proved more, namely
that for a fixed value of $k\ge 1$ the number of positive integers
satisfying the above property has bounded gaps. Specifically,
there exists a computable constant $C(k)$, depending only on $k$, such
that every
interval of length $C(k)$ of positive numbers contains a positive integer
$n$ satisfying the above property. This result has been extended in \cite{CZ}
to the following setting: There exists a positive constant $C(k)$ depending
only on $k$, such that if $\varepsilon_i\in \{0,1\}$ (for $i=1,\dots,k$)
are such that there exists at least one positive integer $n$ satisfying the
congruences $e_{p_i}(n)\equiv \varepsilon_i\pmod 2$ for
all $i=1,\dots,k$,  then there exist infinitely many such positive integers
$n$, and the set of such positive integers has gaps bounded by $C(k)$.
This extends the result of Berend on the original problem
addressed by Erd\H os and Graham, because there $n=1$ clearly
satisfies the above congruences with $\varepsilon_i=0$ for all
$i=1,\dots, k$. The authors of \cite{CZ} mention that for $2\le k\le 5$
every pattern $\varepsilon_i\in \{0,1\}$ for $1\le i\le k$ appears
at least once (hence, infinitely often by the above result), and they
conjectured that every pattern of length $k$ appears for all values of
$k\ge 1$. This conjecture was recently proved by Chen \cite{Chen}.

\par\noindent In \cite{S}, a stronger conjecture is stated, and this conjecture
was mentioned by Yong-Gao Chen in his talk at the ICM2002 in Beijing.
\begin{conj}[Sander \cite{S}]
Let $p_1,\dots,p_k$ be distinct primes, and let
$\varepsilon_1,\dots,\varepsilon_k\in \{0,1\}$. Then
$$|\{0\le n<N: e_{p_i}(n)\equiv \varepsilon_i\pmod 2,~1\le i\le k\}|
\sim {\frac{N}{2^k}} \quad {\rm as}\ N\rightarrow \infty.$$
\end{conj}

\par\noindent In \cite{S}, it is shown that the above conjecture holds for
$k=1$, and it is also shown that every pattern of length $2$ appears
with two arbitrary primes $p_1$ and $p_2$.

\par\noindent In this paper, we first address the following more general
conjecture:
\begin{conj}
Let $p_1,\dots,p_k$ be distinct primes, $m_1,\dots,m_k$
be arbitrary positive integers $(\ge 2)$, and $0\le a_i\le m_i-1$ for $i=1,\dots,k$
be arbitrary residue classes modulo $m_i$. Then
$$|\{0\le n<N:~e_{p_i}(n)\equiv a_i\pmod {m_i},~1\le i\le k\}|\sim
{\frac{N}{m_1\dots m_k}}\quad {\rm as}\ N\rightarrow \infty. $$
\end{conj}

\par\noindent Thus, Conjecture 2 says that for fixed prime numbers
$p_1,\dots, p_k$ and for fixed integers $m_1,\dots, m_k$ the numbers
$(e_{p_1}(n),\dots, e_{p_i}(n))$ should be uniformly distributed modulo
$(m_1,\dots, m_k)$, which is just a natural extension of Conjecture~~1.

\vskip 3pt
\par\noindent  We give the following partial
result which resolves the above Conjecture 2 under some technical assumptions.

\vskip 5pt
\par\noindent {\bf Theorem 1.}
{\em Let $p_1,\dots, p_k$ be distinct primes,
$m_1,\dots, m_k$ be arbitrary positive integers $(\ge 2)$ and
$0\le a_i\le m_i-1$ for $i=1,\dots ,k$ be arbitrary residue classes modulo
$m_i$. Assume further that $p_i{\not|}m_i$ for $i=1,\dots,k$. Then
$$|\{0\le n<N: ~e_{p_i}(n)\equiv a_i\pmod{m_i},~1\le i\le k\}|=
{\frac{N}{m_1\dots m_k}}+O(N^{1-\delta}),$$with $\delta:=1/(120k^2p^{3m}m^2)$,
where $m:=\max\{m_i:~1\le i\le k\}$, and $p:={\max}\{p_i:~1\le i\le k\}$.}

\vskip 5pt
\par\noindent Notice that Theorem 1 not only proves Conjecture 2 under the
particular assumptions that $p_i{\not|m_i}$ for $i=1,\dots, k$, but it
even gives an upper bound on the size of the error term. In particular,
Theorem 1 proves Conjecture 1 for all finite sets of prime numbers not
containing the prime number~~$2$.

\par\noindent In the last section of his paper \cite{B}, Berend asks the following
question: Assume that $2=p_1<p_2<\dots$ is the increasing sequence
of all the prime numbers. Given $k\ge 1$ and an arbitrary pattern
$\varepsilon_i\in \{0,1\}$ for $i=1,2,\dots,k$, how far does one need
to go (as a function of $k$) in order to insure that one finds a number
$n$ so that
$e_{p_i}(n)\equiv \varepsilon_i~({\rm mod}~2)$. Since there are exactly
$2^k$ such patterns, we expect such a number to be larger than $2^k$. Put it
differently, find the fastest growing function $f(x)$ such that for large $x$
there exists $n<x$ such that the congruences
$e_{p_i}(n)\equiv \varepsilon_i\pmod 2$
hold for all $i\le f(x)$, where $\varepsilon_i$ is an arbitrary function
defined on the set of positive integers with values in $\{0,1\}$. While
we are unable to prove that $f(x)$ can be taken of the form $c\log x$
with some constant $c$ (the optimal one being $1/\log 2$), we can
give the following lower bound for $f(x)$.

\par\noindent {\bf Theorem 2.}
{\em Let $3=p_1<p_2<\ldots$ be
the sequence of {\it odd} prime numbers and let $(\varepsilon_i)_{i\ge 1}$
be an arbitrary sequence taking only the values $0$ and $1$.
Then there exists an absolute constant
$c_1$ such that for each large
positive real number $x$, there exists a positive integer $n<x$
such that all congruences
$$
e_{p_i}(n)\equiv \varepsilon_i\pmod 2,\ i=1,\ldots,k(x)
$$
hold with $\displaystyle k(x):={\Bigl\lfloor
c_1 \left(\frac{\log x}{\left(\log\log x\right)^6}
\right)^{1/9}}\Bigr\rfloor$.}

\section{ The Proofs}

\par\noindent
In 1999, D.H. Kim (see \cite{K}) published a very important paper
in which he generalized results of Gelfond  \cite{G}, B\'esineau \cite{Be}, and others,
concerning the joint distribution of completely
$q$-additive functions in residue
classes. The results from Kim's paper found immediate applications, an example
of such being a paper by Thuswaldner and Tichy (see \cite{TT}), in which these
two authors use Kim's results
to prove a version of the Erd\H os-Kac theorem for
sets described by congruence properties of systems of
completely $q$-additive functions.


\par\noindent In this paper, we point out that our theorems are immediate
applications of Kim's main result.


\par\noindent We now describe Kim's main result, and explain how it can be
used to prove our results.


\par\noindent Let $q\ge 2$ be any integer. A function $f$ defined on ${\bf N}
\cup \{0\}$ satisfying $f(0)=0$ and $f(aq^k+b)=f(a)+f(b)$ for all integers
$a\ge 1$, $k\ge 1$, and $0\le b<q^k$ is called {\it completely $q$-additive}.
Note that a completely $q$-additive function may be characterized also
as given by the sum of the values of some function, taken over
the base $q$ digits of the argument.
\par\noindent Let ${\bf q}:=(q_1,\dots,q_k)$, and ${\bf m}:=(m_1,\dots,m_k)$
be $k$-tuples of integers satisfying $q_i,~m_i\ge 2$ for $i=1,\dots, k$,
and $\gcd(q_i,q_j)=1$ for $i\ne j$. For each $i$, let $f_i$ be a
completely $q_i$-additive function with integer values, and write
${\bf f}:=(f_1,\dots, f_k)$. Define
$$
\label{(1)}
F_i:=f_i(1),\eqno (1)
$$
$$
\label{(2)}
d_i:=\gcd\left(m_i, (q_i-1)F_i, f_i(r)-rF_i~(2\le r\le q_i-1)\right),\eqno (2)
$$
and let
${\bf F}:=(F_1,\dots, F_k)$ and ${\bf d}:=(d_1,\dots, d_k)$. For brevity,
we write ${\bf f}(n)\equiv {\bf a}\pmod{\bf m}$ if the congruence
$f_i(n)\equiv a_i\pmod{m_i}$ holds for all $i=1,\dots,k$. Assume
further that ${\rm gcd}(F_i,d_i)=1$ for all $i=1,\dots,k$, and that
${\rm gcd}(d_i,d_j)=1$ for all $i\ne j$.
Kim's theorem says the following:

\par\noindent {\bf Theorem K \cite{K}.}
{\em With the previous assumptions and notations, the estimate
$$
\label{(3)}
|\{0\le n<N:~{\bf f}(n)\equiv {\bf a}\pmod{{\bf m}}\}|=
{\frac{N}{m_1\dots m_k}}+O(N^{1-\delta_1})\eqno (3)$$
holds as $N$ goes to
infinity, and for
all $k$-tuples of residue classes ${\bf a}$ modulo ${\bf m}$, where
$\delta_1:=1/(120k^2q^3m^2)$, with $q:={\rm max}\{q_i:~1\le i\le k\}$, and
$m:={\rm max}\{m_i:~1\le i\le k\}.$}

\begin{proof}[Proof of Theorem $1$]
To apply Kim's theorem, let $i$ be any fixed index in
$\{1,\dots, k\}$. Let $\lambda_i$ be the minimal positive integer such that
the congruence
${\displaystyle{{\frac{p_i^{\lambda_i}-1}{p_i-1}}\equiv 0\pmod{m_i}}}$
holds. Such a value $\lambda_i$
exists because $p_i{\not|}m_i$. Clearly, $\lambda_i\ge 2$
because $m_i\ge 2$. To estimate $\lambda_i$ from above,
write $m_i:=m_i'm_i''$, where $m_i'$ and $m_i''$ are coprime,
all the prime factors of $m_i'$ divide $p_i-1$, and $m_i''$ is coprime to
$p_i-1$.
It is clear that $m_i'$ and $m_i''$
are uniquely determined. By the well-known divisibility properties of
Lucas sequences (see \cite{BHV}), it follows easily that the number
$\mu_i:=m_i'\phi(m_i'')$
satisfies the condition that ${\displaystyle{{\frac{p_i^{\mu_i}-1}{p_i-1}}\equiv 0\pmod{m_i}}}$.
Here, for an arbitrary positive integer
$n$ we used $\phi(n)$ to denote the Euler $\phi$ function of $n$.
It is easy to see that $m_i'\phi(m_i'')\ge 2$ holds for all
positive integers $m_i\ge 2$ (if $m_i=2$, then $p_i$ is odd,
in which case
$\lambda_i=m_i'=m_i=2$, while if $m_i>2$, then either $m_i'\ge 2$,
or $\phi(m_i'')=\phi(m_i)\ge 2$). In particular, $\lambda_i\le m_i\le m$, where
the number $m$ is defined in the statement of Theorem 1.

\par\noindent We now write $q_i:=p_i^{\lambda_i}$. Notice that
${\rm gcd}(q_i,q_j)=1$ holds for all $i\ne j$, and that
$q:={\rm max}\{q_i:~1\le i\le k\}\le p^m$, where $p$ and $m$ are defined
in the statement of Theorem 1.

\par\noindent Define the completely
$q_i$-additive function $f_i$ as follows. Let $a$ be an integer
in the interval $0\le a\le q_i-1$, and write it in base $p_i$ as
$$
\label{(4)}
a:=a_0+a_1p_i+\dots+a_{\lambda_i-1}p_i^{\lambda_i-1}\hskip 20pt {\rm with}~
0\le a_j\le p_i-1~{\rm for}~j=0,\dots,\lambda_i-1.\eqno (4)$$
Set
$$
\label{(5)}
f_i(a):=a_0{\frac{p_i^0-1}{p_i-1}}+a_1{\frac{p_i-1}{p_i-1}}+\dots+a_{\lambda_i-1}
{\frac{p_i^{\lambda_i-1}-1}{p_i-1}},\eqno (5)$$
and extend $f_i$ in the obvious
way to all the non-negative integers in such a way that it becomes
completely $q_i$-additive. A compact formula of $f_i$ in all non-negative
integers is obtained as follows. Let $n\ge 0$, and write it in base $p_i$ as
$$
\label{(6)}
n:=n_0+n_1p_i+\dots+n_tp_i^t\hskip 20pt {\rm with}~0\le n_j\le p_i-1~{\rm for~all}~
j=0,\dots, t.\eqno (6)
$$
For every non-negative integer $j$ write ${\overline{j}}\equiv j \pmod{\lambda_i}$ (the least residue). Then
$$
\label{(7)}
f_i(n)=\sum_{j=0}^t n_j {\frac{p_i^{\overline {j}}-1}{p_i-1}}.\eqno (7)
$$
The next observation is that $f_i(n)\equiv e_{p_i}(n)\pmod{m_i}$. Indeed, it
is well-known that, if the base $p_i$ representation of $n$ is given by (6),
then
$$
\label{(8)}
e_{p_i}(n)=\sum_{j=0}^t n_j {\frac{p_i^j-1}{p_i-1}}.\eqno (8)$$
Thus, comparing
(7) with (8), it suffices to show that
$${\frac{p_i^j-1}{p_i-1}}\equiv {\frac{p_i^{\overline {j}}-1}{p_i-1}}\pmod{m_i},$$
which is equivalent to
$${\frac{p_i^{j-{\overline {j}}}-1}{p_i-1}}\equiv 0\pmod{m_i},$$
and the last congruence is obvious by the definition of $\lambda_i$
and by the fact that $j-{\overline {j}}$ is a multiple of $\lambda_i$.

\par\noindent Having concluded that the completely $q_i$-additive function
$f_i$ represents precisely $e_{p_i}$ modulo $m_i$, in order to complete
the proof of Theorem 1, it suffices, via Theorem K, to verify that the
functions $f_i$ satisfy the assumptions of Theorem K. But
it is clear that by choosing $r:=p_i$, then $2\le r\le q_i-1$ (because
$\lambda_i\ge 2$), and with such an $r$ we have
$f_i(r)-rF_i=f_i(p_i)=1$ (notice that $F_i=f_i(1)=0$).
Thus, formula (2) tells us that $d_i=1$, therefore that all the relations
${\rm gcd}(F_i,d_i)=1$ for $i=1,\dots, k$, and ${\rm gcd}(d_i,d_j)=1$ for
$i\ne j$ hold. All the assertions of Theorem 1 can now be read off from
Theorem K.
\end{proof}

\begin{proof}[Proof of Theorem $2$]
Given $k$, we shall apply Kim's theorem with the error term to find
$N:=N(k)$ in such a way as to make sure that
there exists $n<N$ such that $e_{p_i}(n)\equiv \varepsilon_i\pmod 2$. In order
to do so, we shall need first to find the dependence of the constants
understood in $O$ from Theorem K as a function of the data $k,~m,$ and $q$.
A close analysis of Kim's arguments
points out that the error term appearing in Theorem K can be made explicit
by going through the arguments from Proposition 1 and Proposition 2 in \cite{K}.
In both Propositions 1 and 2 in \cite{K}, the constant understood in the error
term can be taken to be of the type $c_2 k^{1/2}q$, where $c_2$ is an
absolute constant. Thus, with Theorem 1, we have that $m:=2$,
$p_k$ is the $k$th odd prime, $q:=p_k^2$, and
 $\delta:=1/(480k^2p_k^6)$. By Theorem 1 and Kim's theorem
with the explicit dependence of $O$ on the initial data, it follows
that there exists an absolute constant $c_3$ such that whenever the inequality
$$
\label{(9)}
{{N}\over {2^k}}>c_3 k^{1/2}p_k^2 N^{1-\delta}\eqno (9)$$
holds, then there must exist a positive integer
$n<N$ such that $e_{p_i}(n)\equiv \varepsilon_i~({\rm mod}~2)$ holds
for all $i=1,2,\dots,k$. Inequality (9) is equivalent to
$$\label{(10)}
\log N>{\frac{1}{\delta}}\cdot \log(c_3 2^k k^{1/2} p_k^2)=480k^2p_k^6
(\log c_3+k\log 2+0.5\log k+2\log p_k).\eqno (10)
$$
With the Prime Number Theorem, we have
$p_k=(1+o(1))k \log k$, and so inequality (10) yields
$$\label{(11)}
\log N>480(\log 2)(1+o(1))k^9\log^6 k.\eqno (11)
$$
Setting $x=N$ in the left hand side
of the above inequality, and expressing $k$
as a function of $x$, we see that there exists indeed an absolute
constant $c_1$ so that by setting $k$ to be the largest integer less than or
equal than ${\displaystyle{{c_1\Bigl({\frac{\log x}{(\log\log x)^6}}\Bigr)^{1/9}}}}$ then inequality (11) holds.
\end{proof}

\section{Comments}

\par\noindent In this note, we just pointed out how problems about the
distributions of exponents of (fixed) primes appearing in the prime
power factorization of $n!$ in residue classes should be tackled via
the general theory of joint distributions of completely $q$-additive
functions in residue classes. Indeed, by the simple observation of
treating this problem in this way, we pointed out that a result far more
general than any other results available in the literature on this
topic can be inferred in a straightforward way from Kim's results. In this
spirit, we assert that it is probably not too hard to prove Conjecture 2
in its full generality. However, in order to do so,
it is probably far more interesting
and worthwhile to try to
prove an extension of Kim's results to the following setting.

\par\noindent Let ${\bf q}=(q_1,\dots,q_k)$, ${\bf m}:=(m_1,\dots,m_k)$, and
${\bf f}=(f_1,\dots,f_k)$ be as in Kim's Theorem.
Let further ${\bf u}=(u_1,\dots,u_k)$ and
${\bf v}=(v_1,\dots,v_k)$, where $u_i\ge 1$ and $v_i$ are non-negative
integers for $i=1,\dots, k$. Write ${\bf f}({\bf u}n+{\bf v})=
(f_1(u_1n+v_1),\dots,f_k(u_kn+v_k))$. We argue that it would be worthwhile
to study the distribution of the positive integers $n$ such that ${\bf f}({\bf
u}n+{\bf v})\equiv {\bf a}\pmod{\bf m}$, and to conclude that, under
certain natural arithmetical conditions on ${\bf q},~{\bf m},~{\bf u},~
{\bf v}$ and ${\bf f}$, the numbers ${\bf f}({\bf u}n+{\bf v})$ are uniformly
distributed in arithmetical progressions. For example, let us see how one
would attempt to include the number $2$ into the picture in order to prove
Conjecture 1 in its full generality.
If $p_i$ is an odd prime, then define the completely
$q_i$-additive function $f_i$ as in the present paper. When $p_i=2$,
then if one writes $n$ binary as $n:=n_0+2n_1+\dots +2^tn_t$, with
$n_i\in \{0,1\}$ for $i=0,\dots, t$, then
$e_{p_i}(n)\equiv n_1+\dots+n_t\pmod 2$.
Fix $\varepsilon_i\in \{0,1\}$, and
assume that $p_1=2$. For $i=1$, define $f_1(n)$ to be the sum of the digits
of $n$ in base $2$. Consider the system of congruences:
$$f_i(2n)\equiv \varepsilon_i\hskip 20pt {\rm for}~i=1,\dots,k,\eqno (12)$$
and
$$f_1(2n+1)\equiv \varepsilon_1+1,\hskip 10pt f_i(2n+1)\equiv
\varepsilon_i\hskip 20pt {\rm for}~i=2,\dots, k.\eqno (13)$$If Kim's Theorem could
be extended as we pointed out above under pertinent assumptions on
${\bf q},~{\bf m},~{\bf u},~{\bf v}$ and ${\bf f}$, and if such
pertinent assumptions were fulfilled for the systems of congruences
(12) and (13), (with the obvious choices on ${\bf q},~{\bf m},~{\bf u},~
{\bf v}$ and ${\bf f}$),
then one would get a positive answer to Conjecture 1.
By a procedure similar to the one indicated above, one could also
get a positive answer to Conjecture 2 via such an extension of Kim's results.

{\bf Acknowledgements.}
We thank Jean-Marie de Koninck and the referee for suggestions
which improved the presentation of this paper. The first author was supported
in part by grants SEP-CONACYT 37259-E and 37260-E.
The second author was partially supported by a research grant from
the School of Sciences at his institution.

\end{document}